\documentclass[12pt]{amsart}

\usepackage[dvips]{graphicx}

\newtheorem{theorem}{Theorem}[section]
\newtheorem{lemma}[theorem]{Lemma}

\theoremstyle{definition}
\theoremstyle{proposition}

\newtheorem{proposition}[theorem]{Proposition}

\theoremstyle{remark}
\newtheorem{remark}[theorem]{Remark}

\begin{document}

\title[On the additivity of the Thurston--Bennequin invariant]{On the additivity of the Thurston--Bennequin
invariant  of Legendrian knots}

\author{Ichiro TORISU}
\address{Department of Computer Science and Engineering, Faculty of Engineering and 
Resource Science, Akita University, Tegata, Akita 010-8502, Japan }
\email{torisu@math.akita-u.ac.jp}



\date{February 26, 2001}

\begin{abstract}
In this article,  we consider the maximal value of 
the Thurston--Bennequin invariant of Legendrian knots which topologically represent a 
fixed knot type in the
standard contact 3-space  and  we prove a formula of the value  under the connected sum
operation of knots.
\end{abstract}

\maketitle

\section{Introduction}

The {\it standard contact structure} $\xi_0$ on 3-space $\Bbb R^3=\{(x,y,z)\}$ is the plane field  on $\Bbb
R^3$   given by  the kernel of the 1-form $dz-ydx$. A {\it Legendrian knot} $K$ in 
the contact manifold $(\Bbb R^3,\xi_0)$ is a knot which is
everywhere tangent to the  contact structure $\xi_0$.
The {\it Thurston--Bennequin invariant} $tb(K)$ of a Legendrian knot $K$  in $(\Bbb R^3,\xi_0)$ 
is the linking number
of
$K$ and  a knot $K'$ which is obtained by moving $K$ slightly along the vector field $\frac{\partial}{\partial
z}$.  For a topological knot type $k$ in $\Bbb
R^3$, the {\it maximal Thurston--Bennequin invariant} $mtb(k)$ is defined
to  be the maximal value of $tb(K)$, where $K$ is a Legendrian knot which topologically represents $k$. 
For any $k$, by the {\it Bennequin's inequality} in [1], we know that $mtb(k)$ is an integer (i.e. not 
$\infty$). There are several computations of $mtb(k)$ (for example, see [3], [5], [8], [9], [10]). 

In this paper, we prove the following theorem. 

\begin{theorem}
Let $k_1\sharp k_2$ be the connected sum of topological knots $k_1$ and $k_2$ in $\Bbb R^3$. 
Then $mtb(k_1\sharp k_2)=mtb(k_1)+mtb(k_2)+1$.

\end{theorem}

\begin{remark}
After writing this paper, the author was informed that J. Etnyre and K. Honda [4] have also obtained 
a result on connected sum of Legendrian knots which extensively includes Theorem 1.1. 

\end{remark}

\section{Fronts}

Let $K$ be a Legendrian knot in $(\Bbb R^3,\xi_0=ker(dz-ydx))$. Then a diagram (i.e. projection) of $K$ in $xz$-plane
is called  {\it front} as in Figure 1.

\begin{figure}
\begin{center}
\includegraphics[width=8cm,height=5cm,keepaspectratio,clip]{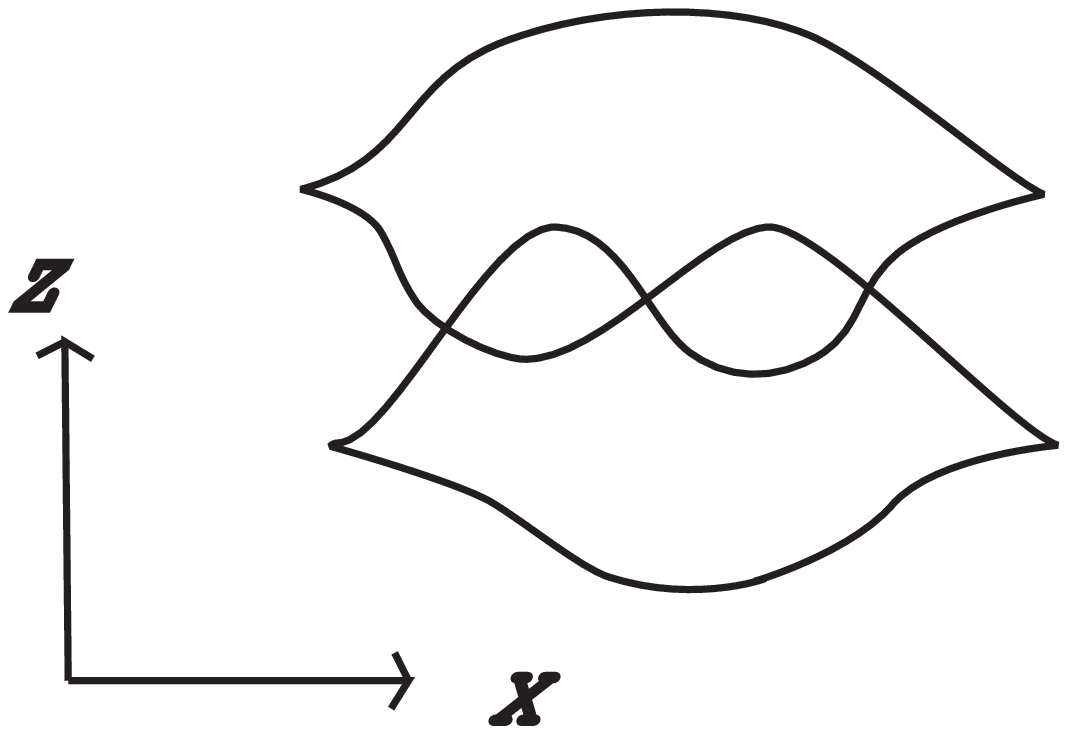}
\end{center}
\caption{}
\label{Figure 1}
\end{figure}%

A front does not have vertical tangents; generically, its only singularities are transverse double points 
and semicubical cusps.  Note that the number of the cusps is even. 
Since $y=\frac{\partial z}{\partial x}$ along $K$, the missing $y$ coordinate is the slope of 
the front. Therefore the front of $K$ is free from selftangencies, and, at a double point, 
the branch with a greater slope is higher along the $y$ axis.
Conversely such a diagram uniquely determines $K$ as its front. So,  
as usual in knot theory, we identify a Legendrian knot $K$ with its front, also denoted by $K$. 

The Thurston--Bennequin invariant $tb(K)$ is computed in terms of the double points and cusps 
of its front. See Figure 2, where $K$ is oriented and the choice of the orientaion is irrelevant for the value 
of $tb(K)$.

\begin{figure}
\begin{center}
\includegraphics[width=12cm,height=5cm,keepaspectratio,clip]{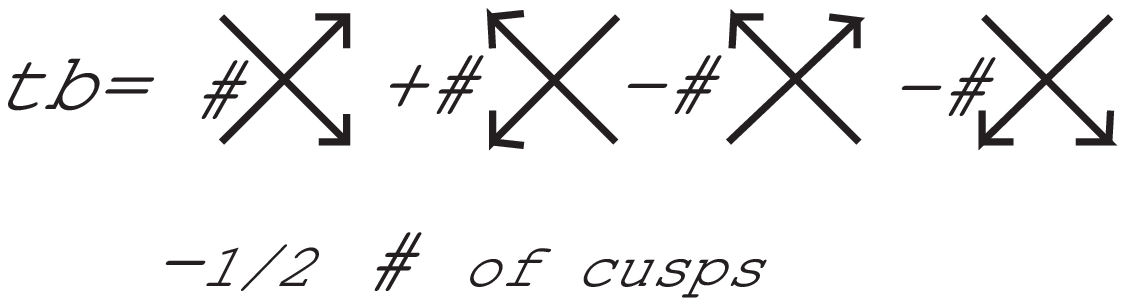}
\end{center}
\caption{}
\label{Figure 2}
\end{figure}%

For example, $tb(K)=-5$ for the front in Figure 1.

\begin{proposition}
For two topological knots $k_1$ and $k_2$, \newline we have $mtb(k_1\sharp k_2)\geq mtb(k_1)+mtb(k_2)+1$. 

\end{proposition}

\begin{proof}

Let $K_1$ and $K_2$ be Legendrian knots whose topological types are $k_1$ and $k_2$, respectively and 
$mtb(k_1)=tb(K_1)$ and $mtb(k_2)=tb(K_2)$.  We also regard $K_1$ and $K_2$ as fronts. 
Further we can assume that $K_1\cap K_2=\emptyset$ and $K_1$ (resp. $K_2$) lies in the left (resp. right) 
region of
$xz$-plane, i.e.  
$\{(x,z)|x<0\}$ (resp. $\{(x,z)|x>0\}$) as in Figure 3. 

\begin{figure}
\begin{center}
\includegraphics[width=12cm,height=7cm,keepaspectratio,clip]{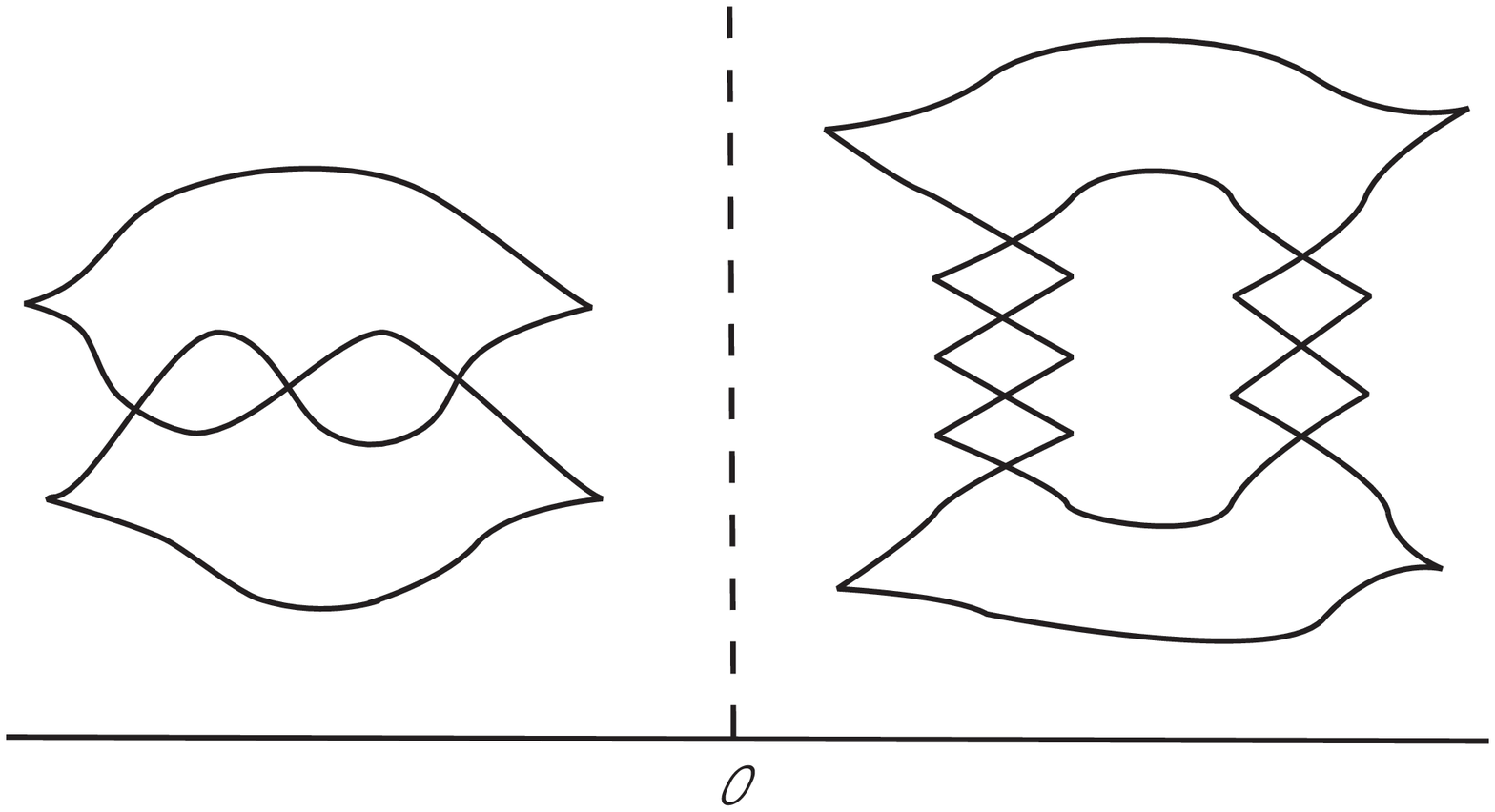}
\end{center}
\caption{}
\label{Figure 3}
\end{figure}%

Then we connect $K_1$ and $K_2$ by joining a right cusp of $K_1$ and a left cusp of $K_2$ as in Figure 4. 

\begin{figure}
\begin{center}
\includegraphics[width=12cm,height=7cm,keepaspectratio,clip]{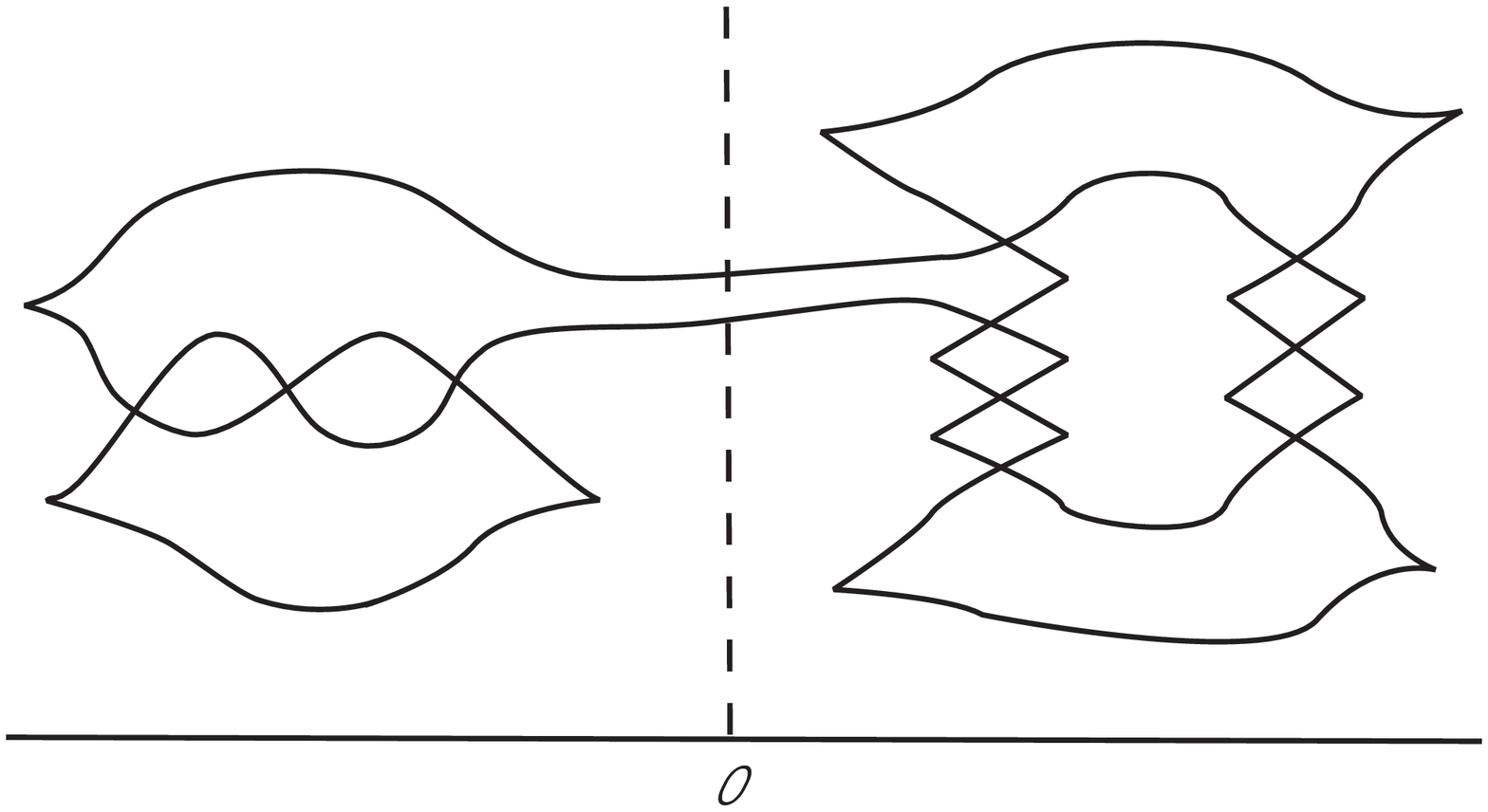}
\end{center}
\caption{}
\label{Figure 4}
\end{figure}%

This procedure produces 
a Legendrian knot whose topological type is $k_1\sharp k_2$ and Thurston--Bennequin invariant is 
$mtb(k_1)+mtb(k_2)+1$. 

This completes the proof.

\end{proof}

\section{Preliminaries from contact topology} 
In this section, 
we recall some basic notions and theorems from recent 3-dimensional contact topology. 
Further, we may assume the reader is familiar with convex surface theory started by E. Giroux in [6]. 
For details
and proofs, see  [2], [3], [6], [7], [8].
Let  $\xi_n=ker(sin(2\pi nz)dx+cos(2\pi nz)dy)$ be the contact structure on a solid torus $V=\{(x,y,z)\in\Bbb
R^{3}_z|x^2+y^2\leq\epsilon\}$, where $n\in\Bbb Z-\{0\}$ and $\Bbb R^{3}_z$ is $\Bbb R^3$ modulo $z\mapsto z+1$. 
The {\it characteristic foliation} on an embedded surface in a contact 3-manifold is   
the singular foliation defined by the intersection of the contact structure and the surface. 
The set of tangents of $\xi_n$ to $\partial V$ forms a disjoint union of two simple closed curves 
on $\partial V$, which are called {\it Legendrian divides}. 
Legendrian divides are leaves of the characteristic foliation on $\partial V$.

The next lemma is proved by a standard Darboux-type argument.
 
\begin{lemma}
 For any Legendrian knot $K$ in $(\Bbb R^3,\xi_0)$, there exists a sufficiently small  
neighborhood $N(K)$ such that 
$(N(K),K,\xi_0)$ is isomorphic to $(V,\{(0,0,z)\},\xi_n)$ for some $n$.
\end{lemma}

Note that in Lemma 3.1, if $K$ is topologically trivial, then $n=tb(K)$.

As $\partial V$ is a {\it convex surface} (i.e. has a contact vector field transverse to $\partial V$),  
the following lemma can be proved by convex surface theory.

\begin{lemma}
Let $T$ be any embedded torus in $(\Bbb R^3,\xi_0)$ and $W$ a solid torus bounded by $T$.   
Suppose the characteristic foliation on $T$ is diffeomorphic to that on $\partial V$ and identifying these, 
the Legendrian divides on $T$ are isotopic to the core curve of $W$ through an isotopy in $W$. 
Then $(W,\xi_0)$ is isomorphic to $(V,\xi_n)$ for some $n$.

\end{lemma}

The following theorem on the classification of topologically trivial Legendrian knots  
due to Y. Eliashberg--M. Fraser [2] is also needed for the proof of Theorem 1.1.

\begin{theorem}
Any topologically trivial Legendrian knot is Legendrian isotopic to one of standard forms 
expressed as fronts in Figure 5.

\end{theorem}

\begin{figure}
\begin{center}
\includegraphics[width=12cm,height=7cm,keepaspectratio,clip]{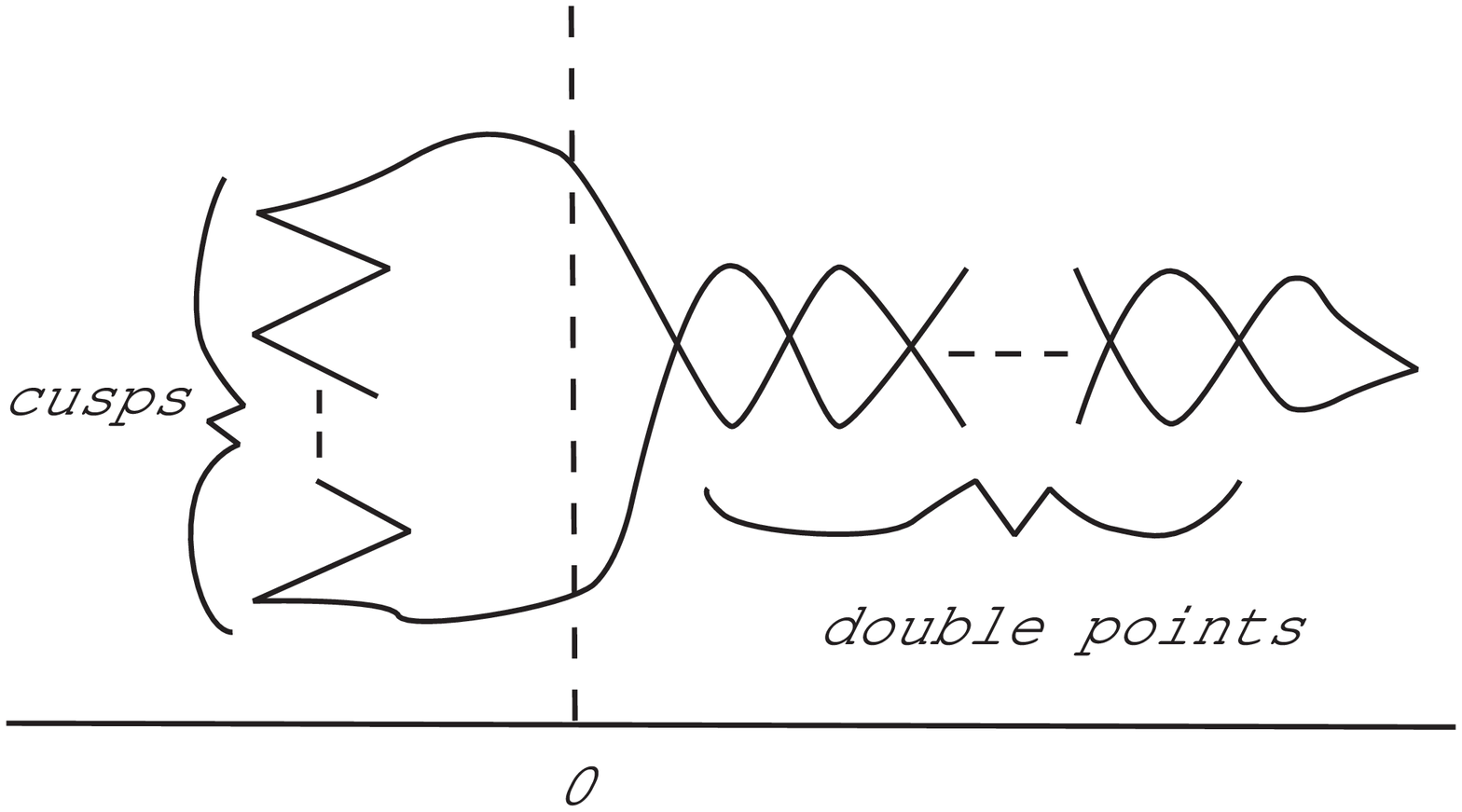}
\end{center}
\caption{}
\label{Figure 5}
\end{figure}%

\section{Proof of Theorem 1.1}

By Propositon 2.1, it is sufficient to show the converse inequality. 

Suppose $\hat K$ is a Legendrian knot in $(\Bbb R^3,\xi_0)$ whose topological type is the connected 
sum of $k_1$ and $k_2$ and its Thurston--Bennequin invariant is maximal. 
By Lemma 3.1, there exists a neighbourhood $N(\hat K)$ of $\hat K$ such that $(N(\hat K),\xi_0)$
 is isomorphc to $(V,\xi_n)$ for some $n$.
Let $B_1$ and $B_2$ be  3-balls in $\Bbb R^3$ such that $B_1$ (resp. $B_2$) splits $\hat K$ into 
the component corresponding to $k_1$ (resp. $k_2$) and $B_1\cap B_2=\emptyset$ (Figure 6).

\begin{figure}
\begin{center}
\includegraphics[width=12cm,height=7cm,keepaspectratio,clip]{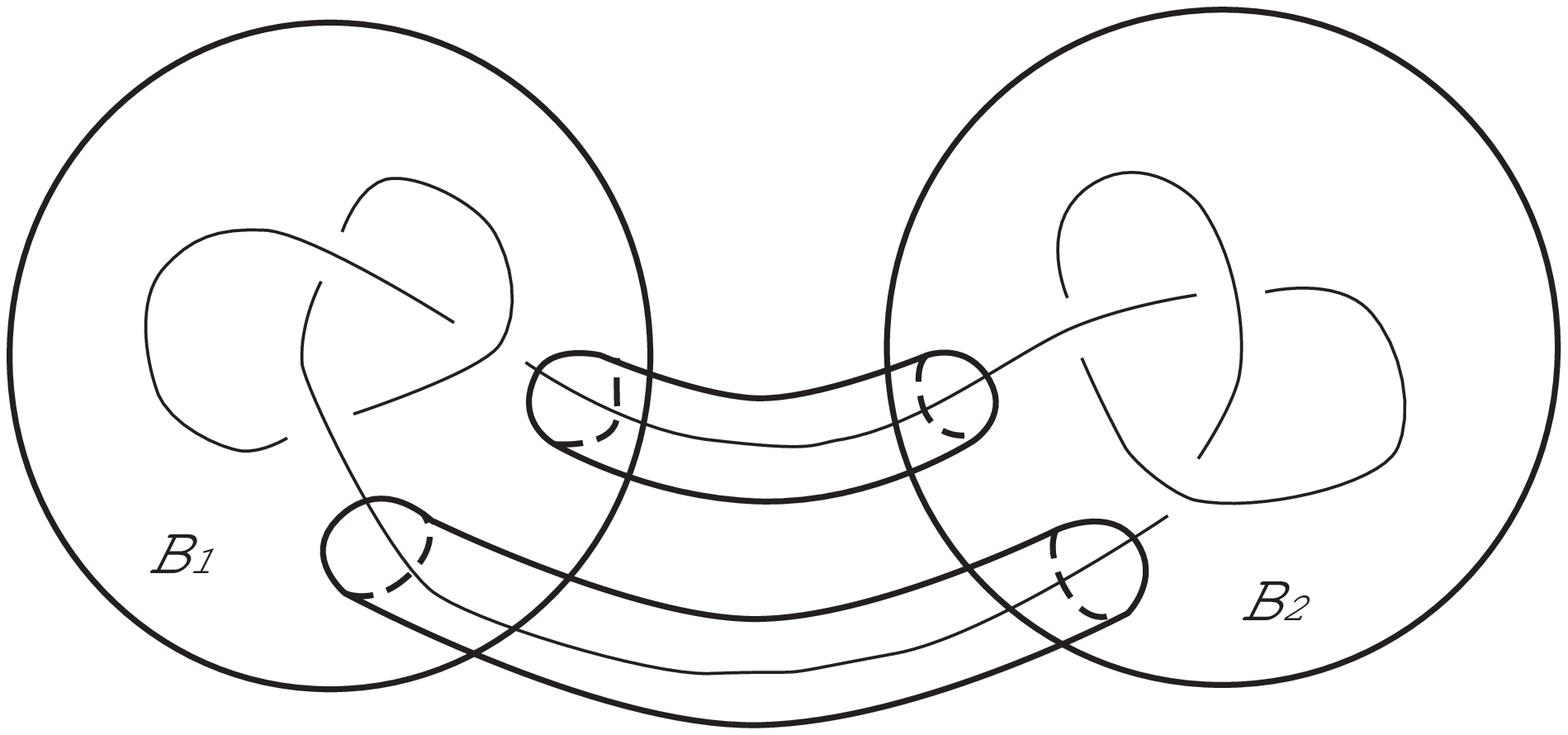}
\end{center}
\caption{}
\label{Figure 6}
\end{figure}%

Further, by convex surface theory, we can assume that (i) $\partial B_1$ and $\partial B_2$ 
 are convex and (ii) $\partial B_1\cap \partial N(\hat K)$ and $\partial B_2\cap \partial N(\hat K)$
 are Legendrian
knots on
$\partial B_1$  and $\partial B_2$, respectively and (iii) each {\it dividing set} on $\partial B_i$
 (i.e. the subset of $\partial B_i$
 consisting of 
tangents  of $\xi_0$ 
 and a contact vector field defining the convex surface) intersects 
$\partial B_i\cap N(\hat K)$ as a diameter of the
disk. 

Then by Edge-Rounding Lemma due to K. Honda in [7], we have a solid torus $W$ such that (i) 
 $W$ equals $B_1\cup B_2\cup N(\hat K)$ except small neighbourhoods of $\partial B_1\cap \partial N(\hat K)$ and
$\partial B_2\cap \partial N(\hat K)$  and (ii) $\partial W$ is a convex surface     
 whose characteristic 
foliation is diffeomorphic to that of $\partial V$. 
By Lemma 3.2, it follows that $(W,\xi_0)$ is isomophic to 
$(V,\xi_n)$ for some $n$. And notice that $W$ is unknotted in $\Bbb R^3$ and hence 
 the core curve $K$ of $W$ which is Legendrian is also unknotted. Further, by a standard 
argument, we can assume that $K$ agrees with $\hat K$ in the region of $N(\hat K)-(B_1\cup 
 B_2)$.
So by Theorem 3.3, $K$ is Legendrian isotopic to one of standard forms in Figure 5. Therefore 
$W$ is also identified with a small neighbourhood of that of the standard form. Further, by a homogeneous 
property of $V$ and a parallel translation of $W$, 
we can assume that a region of $W$ corresponding to $B_1$ (resp. $B_2$) lies in 
$\{(x,y,z)|x<0\}$ (resp. $\{(x,y,z)|x>0\}$). 
Then, identifying $\hat K$ with its front, 
 we can divide $\hat K$ along a vertical line into Legendrian knots $K_1$ and $K_2$ corresponding to 
 $k_1$ and $k_2$, respectively as the converse procedure in the proof of Proposition 2.1. 

Counting the Thurston-Bennequin invariant of $K_1$ and $K_2$, we have 
$tb(\hat K)=mtb(k_1\sharp k_2)=tb(K_1)+tb(K_2)+1$. 
Therefore $mtb(k_1\sharp k_2)\leq mtb(k_1)+mtb(k_2)+1$.

This completes the proof of the main theorem.

\bibliographystyle{amsplain}

\end{document}